\input amstex
\documentstyle{amsppt}

\def\stydate{November 27, 1997}
\immediate\write16{This is DEGT.DEF by A.Degtyarev <\stydate>}
{\edef\temp{\the\everyjob\immediate\write16{DEGT.DEF: <\stydate>}}
\global\everyjob\expandafter{\temp}}

\chardef\tempcat\catcode`\@\catcode`\@=11


\let\le\leqslant
\def\C{{\Bbb C}}
\def\R{{\Bbb R}}
\def\Z{{\Bbb Z}}

\def\Cp#1{\C{\operator@font p}^{#1}}
\def\Rp#1{\R{\operator@font p}^{#1}}
\def\Hom{\qopname@{Hom}}
\def\Ext{\qopname@{Ext}}
\def\Tors{\qopname@{Tors}}

\def\Im{\qopname@{Im}}          
\def\Re{\qopname@{Re}}     %
\def\Ker{\qopname@{Ker}}
\def\Coker{\qopname@{Coker}}
\def\Int{\qopname@{Int}}
\def\Cl{\qopname@{Cl}}
\def\Fr{\qopname@{Fr}}
\def\Fix{\qopname@{Fix}}
\def\tr{\qopname@{tr}}
\def\inj{\qopname@{in}}
\def\id{\qopname@{id}}
\def\pr{\qopname@{pr}}
\def\rel{\qopname@{rel}}
\def\pt{{\operator@font{pt}}}
\def\const{{\operator@font{const}}}
\def\codim{\qopname@{codim}}
\def\cdim{\qopname@{dim_{\C}}}
\def\rdim{\qopname@{dim_{\R}}}
\def\conj{\qopname@{conj}}
\def\rank{\qopname@{rk}}
\def\sign{\qopname@{sign}}
\def\gcd{\qopname@{g.c.d.}}

\def\set<#1|#2>{\bigl\{#1\bigm|#2\bigr\}}


\def\preprint#1{\hrule height0pt depth0pt\kern-24pt%
  \hbox to\hsize{#1}\kern24pt}
\def\today{\ifcase\month\or January\or February\or March\or
  April\or May\or June\or July\or August\or September\or October\or
  November\or December\fi \space\number\day, \number\year}

\def\n@te#1#2{\leavevmode\vadjust{%
 {\setbox\z@\hbox to\z@{\strut\eightpoint#1}%
  \setbox\z@\hbox{\raise\dp\strutbox\box\z@}\ht\z@=\z@\dp\z@=\z@%
  #2\box\z@}}}
\def\leftnote#1{\n@te{\hss#1\quad}{}}
\def\rightnote#1{\n@te{\quad\kern-\leftskip#1\hss}{\moveright\hsize}}
\def\?{\FN@\qumark}
\def\qumark{\ifx\next"\DN@"##1"{\leftnote{\rm##1}}\else
 \DN@{\leftnote{\rm??}}\fi{\rm??}\next@}

\def\centerpage{\dimen@=6.5truein \advance\dimen@-\hsize\hoffset.5\dimen@}
\ifnum\mag>1000 \centerpage\fi

\def\nologo{\let\logo@\relax}

\expandafter\ifx\csname eat@\endcsname\relax\def\eat@#1{}\fi
\expandafter\ifx\csname operator@font\endcsname\relax
 \def\operator@font{\roman}\fi
\expandafter\ifx\csname eightpoint\endcsname\relax
 \let\eightpoint\small\fi

\catcode`\@\tempcat\let\tempcat\undefined

\def\stydate{May 10, 2002}

\chardef\tempcat\catcode`\@
\ifx\undefined\amstexloaded\input amstex \else\catcode`\@\tempcat\fi
\expandafter\ifx\csname amsppt.sty\endcsname\relax\input amsppt.sty \fi
\let\tempcat\undefined

\immediate\write16{This is LABEL.DEF by A.Degtyarev <\stydate>}
\expandafter\ifx\csname label.def\endcsname\relax\else
  \message{[already loaded]}\endinput\fi
\expandafter\edef\csname label.def\endcsname{%
  \catcode`\noexpand\@\the\catcode`\@\edef\noexpand\styname{LABEL.DEF}%
  \def\expandafter\noexpand\csname label.def\endcsname{\stydate}%
    \toks0{}\toks2{}}
\catcode`\@11
\def\labelmesg@ {LABEL.DEF: }
{\edef\temp{\the\everyjob\W@{\labelmesg@<\stydate>}}
\global\everyjob\expandafter{\temp}}

\def\@car#1#2\@nil{#1}
\def\@cdr#1#2\@nil{#2}
\def\eat@bs{\expandafter\eat@\string}
\def\eat@ii#1#2{}
\def\eat@iii#1#2#3{}
\def\eat@iv#1#2#3#4{}
\def\@DO#1#2\@{\expandafter#1\csname\eat@bs#2\endcsname}
\def\@N#1\@{\csname\eat@bs#1\endcsname}
\def\@Nx{\@DO\noexpand}
\def\@Name#1\@{\if\@undefined#1\@\else\@N#1\@\fi}
\def\@Ndef{\@DO\def}
\def\@Ngdef{\global\@Ndef}
\def\@Nedef{\@DO\edef}
\def\@Nxdef{\global\@Nedef}
\def\@Nlet{\@DO\let}
\def\@undefined#1\@{\@DO\ifx#1\@\relax\true@\else\false@\fi}
\def\@@addto#1#2{{\toks@\expandafter{#1#2}\xdef#1{\the\toks@}}}
\def\@@addparm#1#2{{\toks@\expandafter{#1{##1}#2}%
    \edef#1{\gdef\noexpand#1####1{\the\toks@}}#1}}
\def\make@letter{\edef\t@mpcat{\catcode`\@\the\catcode`\@}\catcode`\@11 }
\def\donext@{\expandafter\egroup\next@}
\def\x@notempty#1{\expandafter\notempty\expandafter{#1}}
\def\lc@def#1#2{\edef#1{#2}%
    \lowercase\expandafter{\expandafter\edef\expandafter#1\expandafter{#1}}}
\newif\iffound@
\def\find@#1\in#2{\found@false
    \DNii@{\ifx\next\@nil\let\next\eat@\else\let\next\nextiv@\fi\next}%
    \edef\nextiii@{#1}\def\nextiv@##1,{%
    \edef\next{##1}\ifx\nextiii@\next\found@true\fi\FN@\nextii@}%
    \expandafter\nextiv@#2,\@nil}
{\let\head\relax\let\specialhead\relax\let\subhead\relax
\let\subsubhead\relax\let\proclaim\relax
\gdef\let@relax{\let\head\relax\let\specialhead\relax\let\subhead\relax
    \let\subsubhead\relax\let\proclaim\relax}}
\newskip\@savsk
\let\@ignorespaces\ignorespaces
\def\@ignorespacesp{\ifhmode
  \ifdim\lastskip>\z@\else\penalty\@M\hskip-1sp%
        \penalty\@M\hskip1sp \fi\fi\@ignorespaces}
\def\ignorespaces{\protect\@ignorespacesp}
\def\@bsphack{\relax\ifmmode\else\@savsk\lastskip
  \ifhmode\edef\@sf{\spacefactor\the\spacefactor}\fi\fi}
\def\@esphack{\relax
  \ifx\penalty@\penalty\else\penalty\@M\fi   
  \ifmmode\else\ifhmode\@sf{}\ifdim\@savsk>\z@\@ignorespacesp\fi\fi\fi}
\let\@frills@\identity@
\let\@txtopt@\identyty@
\newif\if@star
\newif\if@write\@writetrue
\def\@numopt@{\if@star\expandafter\eat@\fi}
\def\checkstar@#1{\DN@{\@writetrue
  \ifx\next*\DN@####1{\@startrue\checkstar@@{#1}}%
      \else\DN@{\@starfalse#1}\fi\next@}\FN@\next@}
\def\checkstar@@#1{\DN@{%
  \ifx\next*\DN@####1{\@writefalse#1}%
      \else\DN@{\@writetrue#1}\fi\next@}\FN@\next@}
\def\checkfrills@#1{\DN@{%
  \ifx\next\nofrills\DN@####1{#1}\def\@frills@####1{####1\nofrills}%
      \else\DN@{#1}\let\@frills@\identity@\fi\next@}\FN@\next@}
\def\checkbrack@#1{\DN@{%
    \ifx\next[\DN@[####1]{\def\@txtopt@########1{####1}#1}%
    \else\DN@{\let\@txtopt@\identity@#1}\fi\next@}\FN@\next@}
\def\check@therstyle#1#2{\bgroup\DN@{#1}\ifx\@txtopt@\identity@\else
        \DNii@##1\@therstyle{}\def\@therstyle{\DN@{#2}\nextii@}%
    \expandafter\expandafter\expandafter\nextii@\@txtopt@\@therstyle.\@therstyle
    \fi\donext@}

\newread\@inputcheck
\def\@input#1{\openin\@inputcheck #1 \ifeof\@inputcheck \W@
  {No file `#1'.}\else\closein\@inputcheck \relax\input #1 \fi}

\def\loadstyle#1{\edef\next{#1}%
    \DN@##1.##2\@nil{\if\notempty{##2}\else\def\next{##1.sty}\fi}%
    \expandafter\next@\next.\@nil\lc@def\next@\next
    \expandafter\ifx\csname\next@\endcsname\relax\input\next\fi}

\let\pagebody@\pagebody
\let\pagetop@\empty
\let\pagebot@\empty
\let\@Xend\empty
\def\pagebody{\pagetop@\pagebody@\pagebot@\@Xend}
\let\@Xclose\empty

\newwrite\@Xmain
\newwrite\@Xsub
\def\W@X{\write\@Xout}
\def\make@Xmain{\global\let\@Xout\@Xmain\global\let\end\endmain@
  \xdef\@Xname{\jobname}\xdef\@inputname{\jobname}}
\begingroup
\catcode`\(\the\catcode`\{\catcode`\{12
\catcode`\)\the\catcode`\}\catcode`\}12
\gdef\W@count#1((\lc@def\@tempa(#1)%
    \def\\##1(\W@X(\global##1\the##1))%
    \edef\@tempa(\W@X(%
        \string\expandafter\gdef\string\csname\space\@tempa\string\endcsname{)%
        \\\pageno\\\cnt@toc\\\cnt@idx\\\cnt@glo\\\footmarkcount@
        \@Xclose\W@X(}))\expandafter)\@tempa)
\endgroup
\def\readaux{\bgroup\checkbrack@\readaux@}
\let\begin\readaux
\def\readaux@{%
    \W@{>>> \labelmesg@ Run this file twice to get x-references right}%
    \global\everypar{}%
    {\def\\{\global\let}%
        \def\/##1##2{\gdef##1{\wrn@command##1##2}}%
        \disablepreambule@cs}%
    \make@Xmain{\make@letter\setboxz@h{\@input{\@txtopt@{\@Xname.aux}}%
            \lc@def\@tempa\jobname\@Name\open@\@tempa\@}}%
  \immediate\openout\@Xout\@Xname.aux%
    \immediate\W@X{\relax}\egroup}
\everypar{\global\everypar{}\readaux}
{\toks@\expandafter{\topmatter}
\global\edef\topmatter{\noexpand\readaux\the\toks@}}
\let\@@end@@\end

\def\@Xclose@{{\def\@Xend{\ifnum\insertpenalties=\z@
        \W@count{close@\@Xname}\closeout\@Xout\fi}%
    \vfill\supereject}}
\def\endmain@{\@Xclose@
    \W@{>>> \labelmesg@ Run this file twice to get x-references right}%
    \@@end@@}
\def\disablepreambule@cs{\\\disablepreambule@cs\relax}

\def\include#1{\bgroup
  \ifx\@Xout\@Xsub\DN@{\errmessage
        {\labelmesg@ Only one level of \string\include\space is supported}}%
    \else\edef\@tempb{#1}\clearpage
      \DN@##1 {\if\notempty{##1}\edef\@tempb{##1}\DN@####1\eat@ {}\fi\next@}%
    \DNii@##1.{\edef\@tempa{##1}\DN@####1\eat@.{}\next@}%
        \expandafter\next@\@tempb\eat@{} \eat@{} %
    \expandafter\nextii@\@tempb.\eat@.%
        \relaxnext@
      \if\x@notempty\@tempa
          \edef\nextii@{\write\@Xmain{%
            \noexpand\string\noexpand\@input{\@tempa.aux}}}\nextii@
        \ifx\undefined\@includelist\found@true\else
                    \find@\@tempa\in\@includelist\fi
            \iffound@\ifx\undefined\@noincllist\found@false\else
                    \find@\@tempb\in\@noincllist\fi\else\found@true\fi
            \iffound@\lc@def\@tempa\@tempa
                \if\@undefined\close@\@tempa\@\else\edef\next@{\@Nx\close@\@tempa\@}\fi
            \else\xdef\@Xname{\@tempa}\xdef\@inputname{\@tempb}%
                \W@count{open@\@Xname}\global\let\@Xout\@Xsub
            \openout\@Xout\@tempa.aux \W@X{\relax}%
            \DN@{\let\end\endinput\@input\@inputname
                    \@Xclose@\make@Xmain}\fi\fi\fi
  \donext@}
\def\includeonly#1{\edef\@includelist{#1}}
\def\noinclude#1{\edef\@noincllist{#1}}

\def\arabicnum#1{\number#1}

\def\Romannum#1{\expandafter\uppercase\expandafter{\romannumeral#1}}
\def\alphnum#1{\ifcase#1\or a\or b\or c\or d\else\@ialph{#1}\fi}
\def\@ialph#1{\ifcase#1\or \or \or \or \or e\or f\or g\or h\or i\or j\or
    k\or l\or m\or n\or o\or p\or q\or r\or s\or t\or u\or v\or w\or x\or y\or
    z\else\fi}
\def\Alphnum#1{\ifcase#1\or A\or B\or C\or D\else\@Ialph{#1}\fi}
\def\@Ialph#1{\ifcase#1\or \or \or \or \or E\or F\or G\or H\or I\or J\or
    K\or L\or M\or N\or O\or P\or Q\or R\or S\or T\or U\or V\or W\or X\or Y\or
    Z\else\fi}

\def\ST@P{step}
\def\ST@LE{style}
\def\N@M{no}
\def\F@NT{font@}
\outer\def\newcounter{\checkbrack@{\expandafter\newcounter@\@txtopt@{{}}}}
{\let\newcount\relax
\gdef\newcounter@#1#2#3{{%
    \toks@@\expandafter{\csname\eat@bs#2\N@M\endcsname}%
    \DN@{\alloc@0\count\countdef\insc@unt}%
    \ifx\@txtopt@\identity@\expandafter\next@\the\toks@@
        \else\if\notempty{#1}\global\@Nlet#2\N@M\@#1\fi\fi
    \@Nxdef\the\eat@bs#2\@{\if\@undefined\the\eat@bs#3\@\else
            \@Nx\the\eat@bs#3\@.\fi\noexpand\arabicnum\the\toks@@}%
  \@Nxdef#2\ST@P\@{}%
  \if\@undefined#3\ST@P\@\else
    \edef\next@{\noexpand\@@addto\@Nx#3\ST@P\@{%
             \global\@Nx#2\N@M\@\z@\@Nx#2\ST@P\@}}\next@\fi
    \expandafter\@@addto\expandafter\@Xclose\expandafter
        {\expandafter\\\the\toks@@}}}}
\outer\def\copycounter#1#2{%
    \@Nxdef#1\N@M\@{\@Nx#2\N@M\@}%
    \@Nxdef#1\ST@P\@{\@Nx#2\ST@P\@}%
    \@Nxdef\the\eat@bs#1\@{\@Nx\the\eat@bs#2\@}}
\outer\def\everystep{\checkstar@\everystep@}
\def\everystep@#1{\if@star\let\next@\gdef\else\let\next@\@@addto\fi
    \@DO\next@#1\ST@P\@}
\def\counterstyle#1{\@Ngdef\the\eat@bs#1\@}
\def\advancecounter#1#2{\@N#1\ST@P\@\global\advance\@N#1\N@M\@#2}
\def\setcounter#1#2{\@N#1\ST@P\@\global\@N#1\N@M\@#2}
\def\counter#1{\refstepcounter#1\printcounter#1}
\def\printcounter#1{\@N\the\eat@bs#1\@}
\def\refcounter#1{\xdef\@lastmark{\printcounter#1}}
\def\stepcounter#1{\advancecounter#1\@ne}
\def\refstepcounter#1{\stepcounter#1\refcounter#1}
\def\savecounter#1{\@Nedef#1@sav\@{\global\@N#1\N@M\@\the\@N#1\N@M\@}}
\def\restorecounter#1{\@Name#1@sav\@}

\def\warning#1#2{\W@{Warning: #1 on input line #2}}
\def\warning@#1{\warning{#1}{\the\inputlineno}}
\def\wrn@@Protect#1#2{\warning@{\string\Protect\string#1\space ignored}}
\def\wrn@@label#1#2{\warning{label `#1' multiply defined}{#2}}
\def\wrn@@ref#1#2{\warning@{label `#1' undefined}}
\def\wrn@@cite#1#2{\warning@{citation `#1' undefined}}
\def\wrn@@command#1#2{\warning@{Preamble command \string#1\space ignored}#2}
\def\wrn@@option#1#2{\warning@{Option \string#1\string#2\space is not supported}}
\def\wrn@@reference#1#2{\W@{Reference `#1' on input line \the\inputlineno}}
\def\wrn@@citation#1#2{\W@{Citation `#1' on input line \the\inputlineno}}
\let\wrn@reference\eat@ii
\let\wrn@citation\eat@ii
\def\nowarning#1{\if\@undefined\wrn@\eat@bs#1\@\wrn@option\nowarning#1\else
        \@Nlet\wrn@\eat@bs#1\@\eat@ii\fi}
\def\printwarning#1{\if\@undefined\wrn@@\eat@bs#1\@\wrn@option\printwarning#1\else
        \@Nlet\wrn@\eat@bs#1\expandafter\@\csname wrn@@\eat@bs#1\endcsname\fi}
\printwarning\Protect
\printwarning\label
\printwarning\ref
\printwarning\cite
\printwarning\command
\printwarning\option

{\catcode`\#=12\gdef\@lH{#}}
\def\@@HREF#1{}
\def\@HREF#1#2{\@@HREF{a #1}{\let\@@HREF\eat@#2}\@@HREF{/a}}
\def\@@Hf#1{file:#1} \let\@Hf\@@Hf
\def\@@Hl#1{\if\notempty{#1}\@lH#1\fi} \let\@Hl\@@Hl
\def\@@Hname#1{\@HREF{name="#1"}{}} \let\@Hname\@@Hname
\def\@@Href#1{\@HREF{href="#1"}} \let\@Href\@@Href
\ifx\undefined\pdfoutput
  \csname newcount\endcsname\pdfoutput
\else
  \def\pdflinkattr{attr{/C [0 0.9 0.9]}}
  \let\pdflinkbegin\empty
  \let\pdflinkend\empty
  \def\@pdfHf#1{file {#1}}
  \def\@pdfHl#1{name {#1}}
  \def\@pdfHname#1{\pdfdest name{#1}xyz\relax}
  \def\@pdfHref#1#2{\pdfstartlink \pdflinkattr goto #1\relax
    \pdflinkbegin#2\pdflinkend\pdfendlink}
  \def\@ifpdf#1#2{\ifnum\pdfoutput>\z@\expandafter#1\else\expandafter#2\fi}
  \def\@Hf{\@ifpdf\@pdfHf\@@Hf}
  \def\@Hl{\@ifpdf\@pdfHl\@@Hl}
  \def\@Hname{\@ifpdf\@pdfHname\@@Hname}
  \def\@Href{\@ifpdf\@pdfHref\@@Href}
\fi
\def\@Hr#1#2{\if\notempty{#1}\@Hf{#1}\fi\@Hl{#2}}
\def\@localHref#1{\@Href{\@Hr{}{#1}}}
\def\@countlast#1{\@N#1last\@}
\def\@@countref#1#2{\global\advance#2\@ne
  \@Nxdef#2last\@{\the#2}\@tocHname{#1\@countlast#2}}
\def\@countref#1{\@DO\@@countref#1@HR\@#1}

\def\Href@@#1{\@N\Href@-#1\@}
\def\Href@#1#2{\@N\Href@-#1\@{\@Hl{@#1-#2}}}
\def\Hname@#1{\@N\Hname@-#1\@}
\def\Hlast@#1{\@N\Hlast@-#1\@}
\def\cntref@#1{\global\@DO\advance\cnt@#1\@\@ne
  \@Nxdef\Hlast@-#1\@{\@DO\the\cnt@#1\@}\Hname@{#1}{@#1-\Hlast@{#1}}}
\def\HyperRefs#1{\global\@Nlet\Hlast@-#1\@\empty
  \global\@Nlet\Hname@-#1\@\@Hname
  \global\@Nlet\Href@-#1\@\@Href}
\def\NoHyperRefs#1{\global\@Nlet\Hlast@-#1\@\empty
  \global\@Nlet\Hname@-#1\@\eat@
  \global\@Nlet\Href@-#1\@\eat@}

\HyperRefs{label}
{\catcode`\-11
\gdef\@labelref#1{\Hname@-label{r@-#1}}
\gdef\@xHref#1{\Href@-label{\@Hl{r@-#1}}}
}
\HyperRefs{toc}
\def\@HR#1{\if\notempty{#1}\string\@HR{\Hlast@{toc}}{#1}\else{}\fi}



\def\bftext{\ifmmode\fam\bffam\else\bf\fi}
\let\@lastmark\empty
\let\@lastlabel\empty
\def\lastmark{\@lastmark}
\let\lastlabel\empty
\let\everylabel\relax
\let\everylabel@\eat@
\let\everyref\relax
\def\newlabel{\bgroup\everylabel\newlabel@}
\def\newlabel@#1#2#3{\if\@undefined\r@-#1\@\else\wrn@label{#1}{#3}\fi
  {\let\protect\noexpand\@Nxdef\r@-#1\@{#2}}\egroup}
\def\w@ref{\bgroup\everyref\w@@ref}
\def\w@@ref#1#2#3#4{%
  \if\@undefined\r@-#1\@{\bftext??}#2{#1}{}\else%
   \@xHref{#1}{\@DO{\expandafter\expandafter#3}\r@-#1\@\@nil}\fi
  #4{#1}{}\egroup}
\def\@@@xref#1{\w@ref{#1}\wrn@ref\@car\wrn@reference}
\def\@xref#1{\rom{\@@@xref{#1}}}
\let\xref\@xref
\def\pageref#1{\w@ref{#1}\wrn@ref\@cdr\wrn@reference}
\def\thepage{\ifnum\pageno<\z@\romannumeral-\pageno\else\number\pageno\fi}
\def\label@{\@bsphack\bgroup\everylabel\label@@}
\def\label@@#1#2{\everylabel@{{#1}{#2}}%
  \@labelref{#2}%
  \let\thepage\relax
  \def\protect{\noexpand\noexpand\noexpand}%
  \edef\@tempa{\edef\noexpand\@lastlabel{#1}%
    \W@X{\string\newlabel{#2}{{\@lastmark}{\thepage}}{\the\inputlineno}}}%
  \expandafter\egroup\@tempa\@esphack}
\def\label#1{\label@{#1}{#1}}
\def\fn@P@{\relaxnext@
    \DN@{\ifx[\next\DN@[####1]{}\else
        \ifx"\next\DN@"####1"{}\else\DN@{}\fi\fi\next@}%
    \FN@\next@}
\def\eat@fn#1{\ifx#1[\expandafter\eat@br\else
  \ifx#1"\expandafter\expandafter\expandafter\eat@qu\fi\fi}
\def\eat@br#1]#2{}
\def\eat@qu#1"#2{}
{\catcode`\~\active\lccode`\~`\@
\lowercase{\global\let\@@P@~\gdef~{\protect\@@P@}}}
\def\Protect@@#1{\def#1{\protect#1}}
\def\disable@special{\let\W@X@\eat@iii\let\label\eat@
    \def\footnotemark{\protect\fn@P@}%
  \let\footnotetext\eat@fn\let\footnote\eat@fn
    \let\refcounter\eat@\let\savecounter\eat@\let\restorecounter\eat@
    \let\advancecounter\eat@ii\let\setcounter\eat@ii
  \let\ifvmode\iffalse\Protect@@\@@@xref\Protect@@\pageref\Protect@@\nofrills
    \Protect@@\\\Protect@@~}
\let\notoctext\identity@
\def\W@X@#1#2#3{\@bsphack{\disable@special\let\notoctext\eat@
    \def\chapter{\protect\chapter@toc}\let\thepage\relax
    \def\protect{\noexpand\noexpand\noexpand}#1%
  \edef\next@{\if\@undefined#2\@\else\write#2{#3}\fi}\expandafter}\next@
    \@esphack}
\newcount\cnt@toc
\def\writeauxline#1#2#3{\W@X@{\cntref@{toc}\let\tocref\@HR}
  \@Xout{\string\@Xline{#1}{#2}{#3}{\thepage}}}
{\let\newwrite\relax
\gdef\@openin#1{\make@letter\@input{\jobname.#1}\t@mpcat}
\gdef\@openout#1{\global\expandafter\newwrite\csname tf@-#1\endcsname
   \immediate\openout\@N\tf@-#1\@\jobname.#1\relax}}
\def\@@openout#1{\@openout{#1}%
  \@@addto\readaux@{\immediate\closeout\@N\tf@-#1\@}}
\def\auxlinedef#1{\@Ndef\do@-#1\@}
\def\@Xline#1{\if\@undefined\do@-#1\@\expandafter\eat@iii\else
    \@DO\expandafter\do@-#1\@\fi}
\def\beginW@{\bgroup\def\do##1{\catcode`##112 }\dospecials\do\@\do\"
    \catcode`\{\@ne\catcode`\}\tw@\immediate\write\@N}
\def\endW@toc#1#2#3{{\string\tocline{#1}{#2\string\page{#3}}}\egroup}
\def\do@tocline#1{%
    \if\@undefined\tf@-#1\@\expandafter\eat@iii\else
        \beginW@\tf@-#1\@\expandafter\endW@toc\fi
}
\auxlinedef{toc}{\do@tocline{toc}}

\let\protect\empty
\def\Protect#1{\if\@undefined#1@P@\@\PROTECT#1\else\wrn@Protect#1\empty\fi}
\def\PROTECT#1{\@Nlet#1@P@\@#1\edef#1{\noexpand\protect\@Nx#1@P@\@}}
\def\pdef#1{\edef#1{\noexpand\protect\@Nx#1@P@\@}\@Ndef#1@P@\@}

\Protect\operatorname
\Protect\operatornamewithlimits
\Protect\qopname@
\Protect\qopnamewl@
\Protect\text
\Protect\topsmash
\Protect\botsmash
\Protect\smash
\Protect\widetilde
\Protect\widehat
\Protect\thetag
\Protect\therosteritem
\Protect\Cal
\Protect\Bbb
\Protect\bold
\Protect\slanted
\Protect\roman
\Protect\italic
\Protect\boldkey
\Protect\boldsymbol
\Protect\frak
\Protect\goth
\Protect\dots
\Protect\cong
\Protect\lbrace \let\{\lbrace
\Protect\rbrace \let\}\rbrace
\let\root@P@@\root \def\root@P@#1{\root@P@@#1\of}
\def\root#1\of{\protect\root@P@{#1}}

\def\frills{\ignorespaces\@txtopt@}
\def\frillsnotempty#1{\x@notempty{\@txtopt@{#1}}}
\def\numberline{\@numopt@}
\newif\if@theorem
\let\@therstyle\eat@
\def\@headtext@#1#2{{\disable@special\let\protect\noexpand
    \def\chapter{\protect\chapter@rh}%
    \edef\next@{\noexpand\@frills@\noexpand#1{#2}}\expandafter}\next@}
\let\AmSrighthead@\rightheadtext
\def\rightheadtext{\checkfrills@{\@headtext@\AmSrighthead@}}
\let\AmSlefthead@\leftheadtext
\def\leftheadtext{\checkfrills@{\@headtext@\AmSlefthead@}}
\def\@head@@#1#2#3#4#5{\@Name\pre\eat@bs#1\@\if@theorem\else
    \@frills@{\csname\expandafter\eat@iv\string#4\endcsname}\relax
        \ifx\protect\empty\@N#1\F@NT\@\fi\fi
    \@N#1\ST@LE\@{\counter#3}{#5}%
  \if@write\writeauxline{toc}{\eat@bs#1}{#2{\counter#3}\@HR{#5}}\fi
    \if@theorem\else\expandafter#4\fi
    \ifx#4\endhead\ifx\@txtopt@\identity@\else
        \headmark{\@N#1\ST@LE\@{\counter#3}{\frills\empty}}\fi\fi
    \@Name\post\eat@bs#1\@\ignorespaces}
\ifx\undefined\endhead\Invalid@\endhead\fi
\def\@head@#1{\checkstar@{\checkfrills@{\checkbrack@{\@head@@#1}}}}
\def\@thm@@#1#2#3{\@Name\pre\eat@bs#1\@
    \@frills@{\csname\expandafter\eat@iv\string#3\endcsname}
    {\@theoremtrue\check@therstyle{\@N#1\ST@LE\@}\frills
            {\counter#2}\@theoremfalse}%
    \@DO\envir@stack\end\eat@bs#1\@
    \@N#1\F@NT\@\@Name\post\eat@bs#1\@\ignorespaces}
\def\@thm@#1{\checkstar@{\checkfrills@{\checkbrack@{\@thm@@#1}}}}
\def\@capt@@#1#2#3#4#5\endcaption{\bgroup
    \edef\@tempb{\global\footmarkcount@\the\footmarkcount@
    \global\@N#2\N@M\@\the\@N#2\N@M\@}%
    \def\shortcaption##1{\global\def\sh@rtt@xt####1{##1}}\let\sh@rtt@xt\identity@
    \DN@{#4{\@tempb\@N#1\ST@LE\@{\counter#2}}}%
    \if\notempty{#5}\DNii@{\next@\@N#1\F@NT\@}\else\let\nextii@\next@\fi
    \nextii@#5\endcaption
  \if@write\writeauxline{#3}{\eat@bs#1}{{} \@HR{\@N#1\ST@LE\@{\counter#2}%
    \if\notempty{#5}.\enspace\fi\sh@rtt@xt{#5}}}\fi
  \global\let\sh@rtt@xt\undefined\egroup}
\def\@capt@#1{\checkstar@{\checkfrills@{\checkbrack@{\@capt@@#1}}}}
\let\captiontextfont@\empty

\ifx\undefined\subsubheadfont@\def\subsubheadfont@{\it}\fi
\ifx\undefined\proclaimfont\def\proclaimfont{\sl}\fi
\ifx\undefined\proclaimfont@\let\proclaimfont@\proclaimfont\fi
\def\proclaimfont{\proclaimfont@}
\ifx\undefined\definitionfont@\def\AmSdeffont@{\rm}
    \else\let\AmSdeffont@\definitionfont@\fi
\ifx\undefined\remarkfont@\def\remarkfont@{\rm}\fi

\def\newfont@def#1#2{\if\@undefined#1\F@NT\@
    \@Nxdef#1\F@NT\@{\@Nx.\expandafter\eat@iv\string#2\F@NT\@}\fi}
\def\newhead@#1#2#3#4{{%
    \gdef#1{\@therstyle\@therstyle\@head@{#1#2#3#4}}\newfont@def#1#4%
    \if\@undefined#1\ST@LE\@\@Ngdef#1\ST@LE\@{\headstyle}\fi
    \if\@undefined#2\@\gdef#2{\headtocstyle}\fi
  \@@addto\moretocdefs@{\\#1#1#4}}}
\outer\def\newhead#1{\checkbrack@{\expandafter\newhead@\expandafter
    #1\@txtopt@\headtocstyle}}
\outer\def\newtheorem#1#2#3#4{{%
    \gdef#2{\@thm@{#2#3#4}}\newfont@def#2#4%
    \@Nxdef\end\eat@bs#2\@{\noexpand\revert@envir
        \@Nx\end\eat@bs#2\@\noexpand#4}%
  \if\@undefined#2\ST@LE\@\@Ngdef#2\ST@LE\@{\proclaimstyle{#1}}\fi}}%
\outer\def\newcaption#1#2#3#4#5{{\let#2\relax
  \edef\@tempa{\gdef#2####1\@Nx\end\eat@bs#2\@}%
    \@tempa{\@capt@{#2#3{#4}#5}##1\endcaption}\newfont@def#2\endcaptiontext%
  \if\@undefined#2\ST@LE\@\@Ngdef#2\ST@LE\@{\captionstyle{#1}}\fi
  \@@addto\moretocdefs@{\\#2#2\endcaption}\newtoc{#4}}}
{
\outer\gdef\newtoc#1{{%
    \@DO\ifx\do@-#1\@\relax
    \global\auxlinedef{#1}{\do@tocline{#1}}{}%
    \@@addto\tocsections@{\make@toc{#1}{}}\fi}}}

\toks@\expandafter{\itembox@}
\toks@@{\bgroup\let\therosteritem\identity@\let\rm\empty
  \let\@Href\eat@\let\@Hname\eat@
  \edef\next@{\edef\noexpand\@lastmark{\therosteritem@}}\donext@}
\edef\itembox@{\the\toks@@\the\toks@}
\def\firstitem@false{\let\iffirstitem@\iffalse
    \global\let\lastlabel\@lastlabel}

\let\rosteritemrefform\therosteritem
\let\rosteritemrefseparator\empty
\def\rosteritemref#1{\hbox{\rosteritemrefform{\@@@xref{#1}}}}
\def\local#1{\label@\@lastlabel{\lastlabel-i#1}}
\def\loccit#1{\rosteritemref{\lastlabel-i#1}}
\def\xRef@P@{\gdef\lastlabel}
\def\xRef#1{\@xref{#1}\protect\xRef@P@{#1}}

\def\iref@P@{\gdef\lastref}
\def\itemref#1#2{\rosteritemref{#1-i#2}\protect\iref@P@{#1}}
\def\iref#1{\@xref{#1}\rosteritemrefseparator\itemref{#1}}

\def\eqref#1{\thetag{\@@@xref{#1}}}
\def\tagform@#1{\ifmmode\hbox{\rm\else\rom{\fi
        (\ignorespaces#1\unskip)\iftrue}\else}\fi}

\let\AmSfnote@\makefootnote@
\def\makefootnote@#1{\bgroup\let\footmarkform@\identity@
  \edef\next@{\edef\noexpand\@lastmark{#1}}\donext@\AmSfnote@{#1}}

\def\clearpage{\ifnum\insertpenalties>0\line{}\fi\vfill\supereject}

\def\proof{\checkfrills@{\checkbrack@{%
    \check@therstyle{\@frills@{\demo}{\frills{Proof}}{}}
        {\frills{}\envir@stack\endremark\envir@stack\enddemo}%
  \envir@stack\endproof\ignorespaces}}}
\def\everyendproof{\qed}
\def\endproof{\nofrillscheck{\frills@{\everyendproof}\revert@envir\endproof\enddemo}}

\let\AmSref\ref
\let\AmSrefstyle\refstyle
\let\plaincite\cite
\def\citei@#1,{\citeii@#1\eat@,}
\def\citeii@#1\eat@{\w@ref{#1}\wrn@cite\@car\wrn@citation}
\def\mcite@#1;{\plaincite{\citei@#1\eat@,\unskip}\mcite@i}
\def\mcite@i#1;{\DN@{#1}\ifx\next@\endmcite@
  \else, \plaincite{\citei@#1\eat@,\unskip}\expandafter\mcite@i\fi}
\def\endmcite@{\endmcite@}
\def\cite#1{\mcite@#1;\endmcite@;}
\PROTECT\cite
\def\refstyle#1{\AmSrefstyle{#1}\uppercase{%
    \ifx#1A\relax \def\@ref@##1{\AmSref\xdef\@lastmark{##1}\key##1}%
    \else\ifx#1C\relax \def\@ref@##1{\AmSref\no\counter\refno}%
        \else\def\@ref@{\AmSref}\fi\fi}}
\refstyle A
\newcounter\refno\null
\newif\ifRefs
\gdef\Refs{\checkstar@{\checkbrack@{\csname AmSRefs\endcsname
  \nofrills{\frills{References}%
  \if@write\writeauxline{toc}{vartocline}{\@HR{\frills{References}}}\fi}%
  \def\ref{\@ref@}\Refstrue\ignorespaces}}}
\let\ref\xref

\newif\iftoc
\pdef\tocbreak{\iftoc\hfil\break\fi}
\def\tocsections@{\make@toc{toc}{}}
\let\moretocdefs@\empty
\def\newtocline@#1#2#3{%
  \edef#1{\def\@Nx#2line\@####1{\@Nx.\expandafter\eat@iv
        \string#3\@####1\noexpand#3}}%
  \@Nedef\no\eat@bs#1\@{\let\@Nx#2line\@\noexpand\eat@}%
    \@N\no\eat@bs#1\@}
\def\MakeToc#1{\@@openout{#1}}
\def\newtocline#1#2#3{\Err@{\Invalid@@\string\newtocline}}
\def\make@toc#1#2{\penaltyandskip@{-200}\aboveheadskip
    \if\notempty{#2}
        \centerline{\headfont@\ignorespaces#2\unskip}\nobreak
    \vskip\belowheadskip \fi
    \@openin{#1}\relax
    \vskip\z@}
\def\contents{\readaux\checkfrills@{\checkbrack@{\@contents@}}}
\def\@contents@{\toc@{\frills{Contents}}\envir@stack\endcontents%
    \def\nopagenumbers{\let\page\eat@}\let\newtocline\newtocline@\toctrue
  \def\@HR{\Href@{toc}}%
  \def\tocline##1{\csname##1line\endcsname}
  \edef\caption##1\endcaption{\expandafter\noexpand
    \csname head\endcsname##1\noexpand\endhead}%
    \ifmonograph@\def\vartoclineline{\Chapterline}%
        \else\def\vartoclineline##1{\sectionline{{} ##1}}\fi
  \let\\\newtocline@\moretocdefs@
    \ifx\@frills@\identity@\def\\##1##2##3{##1}\moretocdefs@
        \else\let\tocsections@\relax\fi
    \def\\{\unskip\space\ignorespaces}\let\maketoc\make@toc}
\def\endcontents{\tocsections@\vskip-\lastskip\revert@envir\endcontents
    \endtoc}

\if\@undefined\selectf@nt\@\let\selectf@nt\identity@\fi
\def\Err@math#1{\Err@{Use \string#1\space only in text}}
\def\textonlyfont@#1#2{%
    \def#1{\RIfM@\Err@math#1\else\edef\f@ntsh@pe{\string#1}\selectf@nt#2\fi}%
    \PROTECT#1}
\tenpoint

\def\newshapeswitch#1#2{\gdef#1{\selectsh@pe#1#2}\PROTECT#1}
\def\shapeswitch#1#2#3{\@Ngdef#1\string#2\@{#3}}
\shapeswitch\rm\bf\bf
\shapeswitch\rm\tt\tt
\shapeswitch\rm\smc\smc
\newshapeswitch\em\it
\shapeswitch\em\it\rm
\shapeswitch\em\sl\rm
\def\selectsh@pe#1#2{\relax\if\@undefined#1\f@ntsh@pe\@#2\else
    \@N#1\f@ntsh@pe\@\fi}

\def\@itcorr@{\leavevmode
    \edef\prevskip@{\ifdim\lastskip=\z@ \else\hskip\the\lastskip\relax\fi}\unskip
    \edef\prevpenalty@{\ifnum\lastpenalty=\z@ \else
        \penalty\the\lastpenalty\relax\fi}\unpenalty
    \/\prevpenalty@\prevskip@}
\def\rom@P@#1{\@itcorr@{\selectsh@pe\rm\rm#1}}
\def\rom{\protect\rom@P@}
\def\Rom@P@#1{\@itcorr@{\rm#1}}
\def\Rom{\protect\Rom@P@}
{\catcode`\-11
\HyperRefs{idx}
\HyperRefs{glo}
\newcount\cnt@idx \global\cnt@idx=10000
\newcount\cnt@glo \global\cnt@glo=10000
\gdef\writeindex#1{\W@X@{\cntref@{idx}}\tf@-idx
 {\string\indexentry{#1}{\Hlast@{idx}}{\thepage}}}
\gdef\writeglossary#1{\W@X@{\cntref@{glo}}\tf@-glo
 {\string\glossaryentry{#1}{\Hlast@{glo}}{\thepage}}}
}
\def\emph#1{\@itcorr@\bgroup\em\ignorespaces#1\unskip\egroup
  \DN@{\DN@{}\ifx\next.\else\ifx\next,\else\DN@{\/}\fi\fi\next@}\FN@\next@}
\def\makequoteactive{\catcode`\"\active}
{\makequoteactive\gdef"{\FN@\quote@}
\gdef\quote@{\ifx"\next\DN@"##1""{\quoteii{##1}}\else\DN@##1"{\quotei{##1}}\fi\next@}}
\let\quotei\eat@
\let\quoteii\eat@
\def\MakeIndex{\@openout{idx}}
\def\MakeGlossary{\@openout{glo}}

\def\endofpar#1{\ifmmode\ifinner\endofpar@{#1}\else\eqno{#1}\fi
    \else\leavevmode\endofpar@{#1}\fi}
\def\endofpar@#1{\unskip\penalty\z@\null\hfil\hbox{#1}\hfilneg\penalty\@M}

\newdimen\normalparindent\normalparindent\parindent
\def\firstparindent#1{\everypar\expandafter{\the\everypar
  \global\parindent\normalparindent\global\everypar{}}\parindent#1\relax}

\@@addto\disablepreambule@cs{%
    \\\readaux\relax
    \\\begin\relax
    \\\readaux@\relax
    \\\@openout\eat@
    \\\@@openout\eat@
    \/\Monograph\empty
    \/\MakeIndex\empty
    \/\MakeGlossary\empty
    \/\MakeToc\eat@
    \/\HyperRefs\eat@
    \/\NoHyperRefs\eat@
}

\csname label.def\endcsname


\def\punct#1#2{\if\notempty{#2}#1\fi}
\def\sppunct{\punct{.\enspace}}
\def\varpunct#1#2{\if\frillsnotempty{#2}#1\fi}

\def\headstyle#1#2{\numberline{#1\sppunct{#2}}\ignorespaces#2\unskip}
\def\headtocstyle#1#2{\numberline{#1\punct.{#2}}\space #2}

\def\specialtocstyle#1#2{#2}
\newcounter\section\null
\newcounter\subsection\section
\newcounter\subsubsection\subsection
\newhead\specialsection[\specialtocstyle]\null\endspecialhead
\newhead\section\section\endhead
\newhead\subsection\subsection\endsubhead
\newhead\subsubsection\subsubsection\endsubsubhead
\def\firstappendix{\global\sectionno0 %
  \counterstyle\section{\Alphnum\sectionno}%
    \global\let\firstappendix\empty}

\def\appendixtocstyle#1#2{\space\numberline{Appendix #1\sppunct{#2}}#2}
\newhead\appendix[\appendixtocstyle]\section\endhead

\let\endAmSdef\enddefinition
\def\proclaimstyle#1#2{\numberline{#2\varpunct{.\enspace}{#1}}\frills{#1}}
\copycounter\thm\subsubsection
\theorem\thm\endproclaim
\proposition\thm\endproclaim
\lemma\thm\endproclaim
\corollary\thm\endproclaim
\definition\thm\endAmSdef
\example\thm\endAmSdef

\def\captionstyle#1#2{\frills{#1}\numberline{\varpunct{ }{#1}#2}}
\newcounter\figure\null
\newcounter\table\null
\newcaption{Figure}\figure\figure{lof}\botcaption
\newcaption{Table}\table\table{lot}\topcaption

\copycounter\equation\subsubsection



\copycounter\thm\subsection
\nologo

\def\Per{\operatorname{Per}}
\def\per{\operatorname{per}}
\def\Fix{\operatorname{Fix}}
\def\Pic{\operatorname{Pic}}
\let\Gf\varphi
\def\ix#1{^{\smash{(#1)}}}
\def\ie,{\emph{i.e.,}}
\def\cf.{\emph{cf}@.}
\def\L{\Cal L}

\def\Dg:{\endgraf{\bf Dg:\enspace}\ignorespaces}

\topmatter
\title
Deformation Finiteness for Real Hyperk\"ahler Manifolds
\endtitle

\author
Alex Degtyarev, Ilia Itenberg, Viatcheslav Kharlamov
\endauthor

\abstract
We show that the number of equivariant
deformation classes of real structures in a given deformation
class of compact hyperk\"ahler manifolds is finite.
\endabstract

\thanks
The second and the third authors are supported by
ANR-05-BLAN-0053-01
\endthanks

\address
Bilkent University\endgraf\nobreak
06800 Ankara, Turkey
\endaddress

\email
degt\@fen.bilkent.edu.tr
\endemail

\address
Universit\'e Louis Pasteur et IRMA (CNRS)\endgraf
7 rue Ren\'e Descartes 67084 Strasbourg Cedex, France
\endaddress

\email
itenberg\@math.u-strasbg.fr
\endemail

\address
Universit\'e Louis Pasteur et IRMA (CNRS)\endgraf
7 rue Ren\'e Descartes 67084 Strasbourg Cedex, France
\endaddress

\email
kharlam\@math.u-strasbg.fr
\endemail


\endtopmatter

\vskip-0.2in
\rightline{To Askold Khovansky with our admiration}

\document

\section{Introduction}\label{intro}

An \emph{irreducible holomorphic symplectic manifold} is a
K\"ahlerian compact simply connected complex manifold~$X$ such that
the space $H^0(X,\Omega^2_X)$ is generated by a nowhere degenerate
holomorphic $2$-form $\omega_X$.

In dimension two the only irreducible holomorphic symplectic
manifolds are $K3$-surfaces. Among examples in higher dimensions
are the Hilbert schemes $\operatorname{Hilb}^n(X)$, where $X$ is a
$K3$-surface, and the generalized Kummer varieties $K^n T$, where
$T$ is a complex two-dimensional torus (see \cite{B}).

Sometimes, the irreducible holomorphic symplectic manifolds are
also referred to as compact hyperk\"ahler manifolds. More precisely,
a {\it compact hyperk\"ahler manifold} is an irreducible holomorphic
symplectic manifold with a fixed K\"ahler class $\gamma_X\in
H^2(X;\R)$. (Here and below, by a \emph{K\"ahler class} we mean a
cohomology class represented by the fundamental form of a K\"ahler
metric.)

Recall that a \emph{real structure} on a complex manifold~$X$ is
an anti-holomorphic involution $\conj\:X\to X$. An
\emph{equivariant deformation} is a Kodaira-Spencer family
$p\:\Cal X\to B$ with real structures on both~$\Cal X$ and~$B$
with respect to which the projection~$p$ is equivariant.
Sometimes we use the same term for the pull-back $p^{-1}(B_\R)$ of
the real part of~$B$. Clearly, locally the two notions coincide
and, thus, give rise to the same equivalence relation.

A real structure
on a compact hyperk\"ahler manifold $(X,\gamma_X)$ is a real
structure~$\conj$ satisfying the additional property
$\conj^*\gamma_X=-\gamma_X$. The usual averaging argument shows
that each real irreducible holomorphic symplectic manifold admits
a skew-invariant K\"ahler class, and two real compact hyperk\"ahler manifolds are
equivariantly deformation equivalent if and only if so are the
underlying real irreducible holomorphic symplectic
manifolds.

Our goal is to prove the following theorem.

\theorem\label{main}
The number of equivariant
deformation classes of real structures in a given deformation
class of compact hyperk\"ahler manifolds is finite.
\endtheorem

This result exhibits what we call the \emph{deformation
finiteness} for real hyperk\"ahler manifolds. As is known, a similar
finiteness result holds for curves and surfaces. (Indeed, the only
birational classes of surfaces for which the result is not found
in the literature, either explicitly or implicitly, are elliptic
surfaces and irrational ruled surfaces, but for the latter two
classes the proof is more or less straightforward.) To our
knowledge, almost nothing is known in higher dimensions.

Theorem~\ref{main} is inspired by the following finiteness result for complex
hyperk\"ahler manifolds, see Huybrechts~\cite{H}: \emph{there
exist at most finitely many deformation types of complex
hyperk\"ahler structures on a fixed underlying smooth manifold.}
Moreover, in many respects, out proof of Theorem~\ref{main}
is similar to Huybrechts' proof of
his complex statement.
The crucial points remain the
Koll\'ar-Matsusaka finiteness theorem~\cite{KM}, the
Demailly-Paun characterization of the K\"ahler cone~\cite{DP}, and the Calabi-Yau
families.

Combining Theorem~\ref{main} with Huybrechts' statement
cited above, one can replace
the deformation class of complex manifolds in
Theorem~\ref{main} with the diffeomorphism type of the underlying smooth
manifold. Alternatively, one can consider the manifolds with a
fixed Beauville-Bogomolov form~$q$ and bounded
constant~$\lambda$, see Section~\ref{periods}.

Note that a finiteness statement similar to Theorem~\ref{main}
holds as well for equivariant deformation classes of holomorphic
involutions (instead of real structures); the proof is literally
the same as that of Theorem~\ref{main}.

\subsection*{Acknowledgements}
We are grateful to
the
{\it Max-Planck-Institut f\"ur Mathematik\/}
and to the
{\it Mathematisches Forschungsinstitut Oberwolfach\/} and its RiP program
for their hospitality and excellent working
conditions which helped us to
accomplish an essential part of this work.

\section{Period spaces}\label{periods}

For any irreducible holomorphic symplectic manifold~$X$ there is a
well defined primitive integral quadratic form
$q_X\:H^2(X;\Z)\to\Z$ with the property that for some positive
constant $\lambda\in\R$ (depending on the differential type
of~$X$) the identity $q_X(x)^n=\lambda(x^{2n}\cap[X])$ holds for
any element $x\in H^2(X;\Z)$. (Here, the primitiveness means that
$q_X$ is not a multiple of another integral form.) This quadratic
form, called the \emph{Beauville-Bogomolov form}, has inertia
indexes $(3, b_2(X)-3)$: it is positive definite on the subspace
$(H^{2,0}\oplus H^{0,2})_\R\oplus\gamma_X\R$ and negative definite
on its orthogonal complement. Recall also that
$H^{2,0}\oplus H^{0,2}$ and
$H^{1,1}$ are orthogonal to each other with respect
to the Beauville-Bogomolov form.

From now on, we fix a deformation
class~${\Cal D}$ of irreducible holomorphic
symplectic manifolds. (Note that
the property of being an irreducible holomorphic
symplectic manifold is stable under K\"ahler deformations,
see~\cite{B},
and that all sufficiently small deformations
of a K\"ahler manifold are K\"ahler.)
This defines the isomorphism type of the
cohomology ring $H^*(X;\Z)$ and the Beauville-Bogomolov form
$(H^2(X;\Z),q_X)$. In particular, this fixes the
constant~$\lambda$ in the above description of the form.

In addition, we fix an abstract integral lattice $(L,q)$ isometric
to $(H^2(X;\Z),q_X)$.
For a particular manifold~$X \in \Cal D$, a choice of an isometry
$\Gf\:(H^2(X;\Z),q_X)\to (L,q)$ is called a \emph{marking} of~$X$.

Denote by $\Per$ the \emph{period domain}
$$
\Per=\big\{z\in L\otimes\C\bigm|
 \text{$q(z)=0$, $q(z+\bar z)>0$}\bigr\}/\C^*\subset
 \Bbb P(L\otimes\C).
$$
The \emph{period map}, denoted by~$\per$, sends a
marked
manifold
$(X,\Gf)$, $X \in \Cal D$,
to the point
$\per(X,\Gf)=\Gf(H^{2,0}(X))\bmod\C^*\in\Per$.
The period
map is known to be surjective (see~\cite{H}), and the following local
Torelli theorem holds (see~\cite{Bo}).

\theorem\label{Torelli-C}
Given
a manifold~$X \in \Cal D$,
there exists a universal \rom(in the class of
local
deformations of~$X$\rom)
local deformation
$$
p\:(\Cal X,\Cal X_0=X)\to(\Cal B,b_0).
$$
Furthermore, for any marking~$\Gf$ of~$X$, the period map
$B\to\Per$, $b\mapsto\per p^{-1}(b)$ \rom(which is well defined
due to the stability mentioned above\rom), is a diffeomorphism of
a neighborhood of the base point~$b_0$ in~$B$ to a neighborhood of
its image $\per(X,\Gf)$ in $\Per$.
\qed
\endtheorem

Recall that a \emph{K\"ahler-Einstein metric} on a manifold~$X$
with complex structure~$I$ is a Ricci flat $I$-invariant
Riemannian metric~$g$ on~$X$ whose
`hermitization' $h(u,v)=g(u,v)-ig(Iu,v)$
is K\"ahler.
The following
fundamental
statement is a well known corollary of the Calabi-Yau
theorem.

\theorem\label{Calabi-Yau}
Let $(X,\gamma_X)$, $X\in\Cal D$, be a compact hyperk\"ahler
manifold.
Then
$\gamma_X$ is represented by a unique closed form~$\rho_I$
such that the corresponding
metric $g(u,v)=\rho_I(u,Iv)$
is K\"ahler-Einstein with respect to the original complex structure~$I$ and
two additional complex structures~$J$ and~$K$ on~$X$
satisfying the relation
$IJ=-JI=K$.
\qed
\endtheorem

Consider
a manifold $X\in\Cal D$ and complex structures $I$, $J$, $K$
as in Theorem~\ref{Calabi-Yau}.
For any triple
$(a,b,c)\in S^2=\bigl\{\bold x\in\R^3\bigm|\mathopen\|\bold
x\mathclose\|=1\bigr\}$, the operator
$\lambda=aI+bJ+cK$
is also
a complex structure
with respect to
which
$g$ is a K\"ahler-Einstein metric.
Thus,
a K\"ahler-Einstein metric on~$X$ defines a whole $2$-sphere of
complex structures, each complex structure coming with a
distinguished
K\"ahler class.
This sphere
is naturally identified with the unit sphere in
the (maximal) positive
definite (with respect to~$q_X$) subspace $V\subset H^2(X;\R)$
spanned
by the corresponding K\"ahler classes $\gamma_I=\gamma_X$, $\gamma_J$,
$\gamma_K$. Alternatively, $V$ is spanned by
$\gamma_X$, $\Re[\omega_X]$, and~$\Im[\omega_X]$
and, thus, depends on~$X$ and~$\gamma_X$ only.
The corresponding sphere of complex structures
is denoted by $\Cal S(X,\gamma_X)$. Each element
$\gamma\in\Cal S(X,\gamma_X)$ is the K\"ahler class of
a compact hyperk\"ahler manifold (defined by~$\gamma$)
deformation equivalent to $(X,\gamma_X)$.

For any real
compact hyperk\"ahler manifold $(X,\gamma_X,\conj)$, $X\in\Cal D$,
the involution $\conj^*$ induced in the cohomology of~$X$
is a $q_X$-isometry.
Furthermore, $\conj^*$
interchanges $H^{2,0}$ and $H^{0,2}$,
preserves $H^{1,1}$, and multiplies~$\gamma_X$ by~$-1$.
Therefore, with respect to~$q_X$,
the $(+1)$- and $(-1)$-eigenlattices of $\conj^*\:H^2(X;\Z)\to H^2(X;\Z)$
have positive inertia
indexes~$1$
and~$2$, respectively.

The following statement is a straightforward consequence of
Theorem~\ref{Calabi-Yau}.

\corollary\label{real.Calabi-Yau}
Let $(X,\gamma_X,\conj)$, $X \in \Cal D$, be
a real compact hyperk\"ahler  manifold,
and let $\Cal S_{\conj}(X,\gamma_X)\subset\Cal S(X,\gamma_X)$ be
the circle fixed by~$\conj^*$.
Then each compact hyperk\"ahler manifold
defined by a class
$\gamma\in\Cal S_{\conj}(X,\gamma_X)$
is real with respect to~$\conj$.
\qed
\endcorollary

\remark{Remark}
Under the assumptions of Corollary~\ref{real.Calabi-Yau},
the holomorphic form~$\omega_X$ on~$X$ can be chosen so that
$\conj^*\omega_X=\overline\omega_X$. With this choice,
$\Cal S_{\conj}(X,\gamma_X)$ is the unit circle in the plane
spanned by~$\gamma_X$ and~$\Im[\omega_X]$.
\endremark

A \emph{real homological type} is a $q$-isometry $c\:L\to L$ whose
$(+1)$-eigensublattice has
positive inertia index~$1$.
Whenever $c$ is understood, we denote
by $L_\pm$ its $(\pm1)$-eigenlattices.
The involution~$c$
induces the map $\omega\mapsto\overline{c(\omega)}$ on the
period space $\Per$. Denote by $\Per(c)$ its fixed point set.
Clearly, $\Per(c)$ consists of pairs $(\omega_+\R,\omega_-\R)$,
where $\omega_+\in L_+\otimes\R$, $\omega_-\in L_-\otimes\R$, and
$q(\omega_+)=q(\omega_-)>0$.

Fix a real homological type~$c$. A \emph{$c$-marking} of a real
manifold $(X,\conj)$, $X\in\Cal D$, is a marking
$\Gf\:(H^2(X;\Z),q_X)\to(L,q)$ commuting with~$c$, \ie, such that
$c\circ\Gf=\Gf\circ\conj^*$. Clearly, the image of a $c$-marked
manifold $(X,\conj,\Gf)$ under the period map~$\per$ belongs to
$\Per(c)$.

\theorem\label{real.Torelli}
Any $c$-marked real manifold
$(X,\conj,\Gf)$, $X\in\Cal D$, admits an equivariant local
deformation over a base diffeomorphic to a neighborhood of the
image $\per(X,\conj,\Gf)$ in $\Per(c)$.
\endtheorem

\proof
Pick a universal local deformation
$p\:(\Cal X, \Cal X_0)\to(\Cal B,b_0)$ of $\Cal X_0=X$
given by Theorem~\ref{Torelli-C}. Due to
the universality, the real structure $\conj\:X\to X$ extends to a
unique fiber preserving anti-holomorphic map $\tilde c\:\Cal
X\to\Cal X$. On the other hand, one has $H^0(X;\Cal
T_X)=H^0(X;\Omega^1_X)=H^{1,0}(X)=0$,
\ie, $X$ has no infinitesimal automorphisms. Hence, $\tilde c$ is
an involution. The restriction of~$p$ to the pull-back of
$\Per(c)$ under the period map $b\mapsto\per p^{-1}(b)$ is the
desired equivariant deformation.
\endproof

Next statement is a corollary of the Demailly-Paun
description~\cite{DP} of the K\"ahler cone; see~\cite{HErra} for
details.

\theorem\label{KC}
Let $\Cal K_X\subset H^2(X;\R)$ be the K\"ahler
cone of a manifold $X\in\Cal D$, and let $\Gf$ be a marking
of~$X$. If the period $\omega=\per(X,\Gf)$ is generic, then the
image $\Gf(\Cal K_X)$ coincides with the positive cone in the
hyperbolic space $\omega^\perp\cap(L\otimes\R)$.
\qed
\endtheorem

\remark{Remark}
More precisely, `generic' in Theorem~\ref{KC} means the following
(see~\cite{F}): no power of~$\omega_X$ should annihilate an
integral homology class of~$X$ and, in addition, at most one point
from each sphere of the form $\Cal S(Y,\gamma)$, $Y\in \Cal D$,
should be removed. This implies that the conclusion of
Theorem~\ref{KC} still holds for a generic {\bf real} period, \ie,
for any real homological type $c\:L\to L$, the set of periods of
real manifolds with maximal K\"ahler cone is dense in $\Per(c)$.
\endremark

The last ingredient of the proof is the following statement due to
J.~Koll\'ar and T.~Matsusaka, see~\cite{KM}.

\theorem
For any integer $n>0$ and any pair of integers~$a$ and~$b$ there are
universal constants~$m$, $N$ such that for any compact complex
$n$-manifold~$X$ and any ample line bundle~$\L$ on~$X$ satisfying
the relations
$c_1(\L)^n=a$ and $c_1(\L)^{n-1}c_1(X)=b$ the bundle~$\L^m$ is very
ample and $\dim\mathopen|L^m\mathclose|\le N$.
\qed
\endtheorem

\corollary\label{KMCor}
Fix a real homological type~$c$
and a
vector $l\in L_-$ with $q(l)>0$.
Then, the real compact hyperk\"ahler
manifolds $(X,\gamma_X,\conj)$, $X\in\Cal D$,
admitting
a $c$-marking $\Gf$ such that $\Gf(\gamma_X)=l$
constitute finitely
many deformation families.
\endcorollary

\proof
After a small perturbation, one can assume that the Picard group
$\Pic X=\per(X,\Gf)^\perp\cap L$ is isomorphic to~$\Z$. Then
$\gamma_X=c_1(\L)$ for some ample line bundle~$\L$, and
Theorem~\ref{KMCor}
applied to $a=l^n$ and $b=0$
gives constants~$m$, $N$ such that the
bundle~$\L^m$ embeds~$X$ to the projective space
$\mathopen|L^m\mathclose|\spcheck$ of dimension at most~$N$, the
degree of the image being~$l^n$.
Since $X$ is simply connected and the Chern class
$c_1(\L^m)$ is $\conj^*$-skew-invariant, the real structure~$\conj$
lifts to a real structure on~$\L$, which descends to a real
structure on $\mathopen|\L^m\mathclose|\spcheck$
(possibly, nonstandard) with respect to
which the above embedding is equivariant. Thus, the manifolds in
question are realized as real submanifolds of bounded degree
of real projective spaces of bounded dimensions; such manifolds
form finitely many deformation families,
each family consisting of compact hyperk\"ahler manifolds
due to the deformation stability.
\endproof

\section{Proof of Theorem~\ref{main}}\label{pr}

\lemma\label{lem}
Let $M$ be a lattice with at least two positive squares,
and let $v_1,v_2\in M\otimes\R$ be two vectors with
$v_1^2>0$ and $v_2^2>0$.
Then there exists a vector $v\in M\otimes\R$
such that the bilinear form is positive definite on both the plane
generated by~$v$ and~$v_1$ and the plane generated by~$v$ and~$v_2$.
\endlemma

\proof
Consider the subspace~$V$ generated by~$v_1$ and~$v_2$.
If it has two positive squares, take $v\in V$; otherwise, $V^\perp$
has a positive square, and take $v\in V^\perp$.
\endproof

\subsection{Proof of Theorem~\ref{main}}
Up to isomorphism, $(L,q)$ admits at most finitely many involutive
isometries $c\:L\to L$. (Indeed, enumerating involutive isometries
of a lattice reduces to enumerating isomorphism classes of
lattices of bounded determinant, \cf.~\cite{Nikulin}, and their
number is finite, see~\cite{Cassels}.)
Thus, it suffices to show that each
real homological type~$c$ is realized by finitely many deformation
families.

Fix~$c$ and pick a class $l\in L_-$ with $q(l)>0$. We will show
that any $c$-marked real compact hyperk\"ahler manifold
$(X,\gamma_X,\conj,\Gf)$, $X\in\Cal D$, is equivariantly
deformation equivalent to a $c$-marked real compact hyperk\"ahler
manifold $(Y,\gamma_Y,\conj,\psi)$ with $\psi(\gamma_Y)=l$. Due to
Corollary~\ref{KMCor}, manifolds~$Y$ with this property constitute
finitely many deformation families.

Consider a $c$-marked real compact hyperk\"ahler manifold
$(X,\gamma_X,\conj,\Gf)$, $X\in\Cal D$. Let
$\gamma=\Gf(\gamma_X)\in L_-\otimes\R$, and denote by
$\omega_\pm\in L_\pm\otimes\R$ the periods of~$X$. To construct
the desired deformation, we use the following three moves of the
triple $(\omega_+,\omega_-,\gamma)$ in the period space:
\roster
\item\local1
a rotation in the plane
spanned
by~$\omega_-$ and~$\gamma$,
\item\local2
a small perturbation of~$\omega_+$ in~$L_+\otimes\R$
and of~$\omega_-$ and~$\gamma$
in $L_-\otimes\R$, so that $\omega_-$ and~$\gamma$ remain
orthogonal to each other,
\item\local3
a move of~$\gamma$ within the
cone $\Gf(\Cal K_X)\cap(L_-\otimes\R)$, where
$\Cal K_X\subset H^2(X;\R)$ is the K\"ahler cone of~$X$.
\endroster
Moves~\loccit1 and~\loccit2 are followed by equivariant deformations
of $(X,\gamma_X,\conj)$ due to
Corollary~\ref{real.Calabi-Yau} and
Corollary~\ref{real.Torelli}, respectively. For move~\loccit3, if the pair
$(\omega_+,\omega_-)$ is sufficiently generic, Theorem~\ref{KC}
implies that $\gamma$ can vary within the whole positive cone
in the hyperbolic space $\omega_-^\perp\subset L_-\otimes\R$.

It remains to start with
$(\omega_\pm\ix0,\gamma\ix0)=(\omega_\pm,\gamma)$ and construct a
sequence of triples $(\omega_\pm\ix{i},\gamma\ix{i})$,
$0\le i\le k$, each obtained from the previous one by
one of the three moves, so that $\gamma\ix{k}=l$. As explained
above, that would give rise to a sequence of $c$-marked compact
hyperk\"ahler manifolds
$(X\ix{i},\gamma_{X\ix{i}},\conj,\Gf)$ deformation equivalent to each
other, and one can take
$(Y,\gamma_Y,\conj,\psi)=(X\ix{k},\gamma_{X\ix{k}},\conj,\Gf)$.

The desired sequence can be constructed as follows. First, perturb
$(\omega_\pm\ix0,\gamma\ix0)$ to a generic triple
$(\omega_\pm\ix1,\gamma\ix1)$ (move~\loccit2). Apply
Lemma~\ref{lem} to get a vector $v\in L_-\otimes\R$ such that both
the plane spanned by~$v$ and $\omega_-\ix1$ and the plane spanned
by~$v$ and~$l$ are positive definite. Since the periods
$\omega_\pm\ix1$ are generic, one can replace~$\gamma\ix1$ with a
vector $\gamma\ix2$ in the plane spanned by~$v$ and $\omega_-\ix1$
(move~\loccit3). Let $\omega_\pm\ix2=\omega_\pm\ix1$, use
move~\loccit1 to produce a triple $(\omega_\pm\ix3,\gamma\ix3)$
with $\omega_-$ a positive multiple of~$v$, and perturb it to a
generic triple $(\omega_\pm\ix4,\gamma\ix4)$ (move~\loccit2) so
that the plane spanned by $\omega_-\ix4$ and~$l$ is still positive
definite. Finally, replace~$\gamma\ix4$ with a vector~$\gamma\ix5$
in the plane spanned by $\omega_-\ix4$ and~$l$ and such that
$q(\gamma\ix5)=q(l)$ (move~\loccit3), let
$\omega_\pm\ix5=\omega_\pm\ix4$, and use move~\loccit1 to produce
a triple $(\omega_\pm\ix6,\gamma\ix6)$ with $\gamma\ix6=l$.
\qed


\Refs\widestnumber\key{KM}

\ref{Be}
\by A.~Beauville
\paper  Vari\'et\'es K\"ahleriennes dont la premi\`ere classe de Chern est nulle
\jour J. Diff. Geom.
\vol
18
\yr 1983
\pages 755-782
\endref\label{B}

\ref{Bo}
\by F.~A.~Bogomolov
\paper Hamiltonian K\"ahler manifolds
\jour Sov. Math. Dokl.
\vol 19
\yr 1978
\pages 1462--1465
\endref\label{Bo}

\ref{DP}
\by J.-P.~Demailly, M.~Paun
\paper Numerical characterization of the K\"ahler cone of a compact K\"ahler manifold
\jour Ann. of Math.
\vol 159
\yr 2004
\pages 1247-1274
\endref\label{DP}

\ref{C}
\by J.~W.~S.~Cassels
\book Rational Quadratic Forms
\publ Academic Press
\publaddr New York, London
\yr 1978
\endref\label{Cassels}

\ref{F}
\by A.~Fujiki
\paper On the de Rham Cohomology Group of a Compact K\"ahler Symplectic Manifold
\jour Adv. Stud. Pure Math.
\vol 10
\yr 1987
\pages 105--165
\endref\label{F}

\ref{H1}
\by D.~Huybrechts
\paper Erratum: Compact hyperk\"ahler manifolds: basic results
\inbook Preprint math. AG/0106014
\endref\label{HErra}

\ref{H2}
\by D.~Huybrechts
\paper   Finiteness results for compact hyperk\"ahler manifolds
\jour     J. reine angew. Math.
\vol  558
\yr 2003
\pages 15-22
\endref\label{H}

\ref{KM}
\by J.~Koll\'ar, T.~Matsusaka
\paper Riemann-Roch type inequalities
\jour Amer. J. Math.
\vol 105
\yr 1983
\pages 229--252
\endref\label{KM}

\ref{N}
\by V.~V.~Nikulin
\paper Integer quadratic forms and some of their geometrical applications
\jour Izv. Akad. Nauk SSSR, Ser. Mat
\vol 43
\yr 1979
\pages 111--177
\lang Russian
\transl\nofrills English transl. in
\jour Math. USSR--Izv.
\vol 43
\yr 1979
\pages 103--167
\endref\label{Nikulin}

\ref{S}
\by Y.-T. Siu
\paper Every K3 surface is K\"ahler
\jour Invent. Math.
\vol 73
\yr
1983
\pages 130-150
\endref\label{S}


\endRefs
\enddocument